\documentclass[12pt,reqno]{amsart}

\usepackage{geometry}
\usepackage[T1]{fontenc}
\usepackage{lmodern,microtype}
\usepackage{amsmath,amssymb,amsthm,mathtools,mathrsfs}
\usepackage{fancyhdr,hyperref}
\hypersetup{hidelinks}

\numberwithin{equation}{section}
\newtheorem{theorem}{Theorem}[section]
\newtheorem{lemma}[theorem]{Lemma}

\theoremstyle{remark}

\newcommand{\dd}{\mathrm{d}}
\newcommand{\SL}{\operatorname{SL}}
\renewcommand{\Im}{\operatorname{Im}}
\newcommand{\vol}{\operatorname{vol}}

\begin{document}
\title[$L^p$-norms of holomorphic cusp forms]
{$L^p$-norms of holomorphic cusp forms}

\author{Xiaolong Han}
\address{Department of Mathematics, California State University, Northridge, CA 91330, USA}
\email{xiaolong.han@csun.edu}

\author{Bingrong Huang}
\address{Data Science Institute \& School of Mathematics, Shandong University, Jinan, Shandong 250100, China}
\curraddr{Department of Mathematics, Princeton University, Fine Hall, Washington Road, Princeton, NJ 08544, USA}
\email{brhuang@sdu.edu.cn}

\date{\today}
\thanks{B.H. was supported by the Key R\&D Program of Shandong Province, China (No. 2026ZLYS01), the Scientific Research Innovation Capability Support Project for Young Faculty (No. SRICSPYF-ZY2025158), and a Mathematics Special Project of Shandong University (No. TQ062025004).}

\begin{abstract}
We determine the optimal orders of local and global $L^p$-norms of
$L^2$-normalized holomorphic cusp forms of large weight on the modular
surface, for fixed $2<p\leq\infty$. Both orders are attained by peak
forms. For normalized Poincar\'e series of weight $k$ and index $m$,
we obtain an asymptotic formula with an explicit constant, uniformly for
$1\leq m\leq c_0k$ with any fixed $c_0<(4\pi)^{-1}$. We also give
sup norm bounds in a wider range of indices.
\end{abstract}
\maketitle

\section{Introduction}\label{sec:introduction}

The mass distribution of holomorphic cusp forms in the weight aspect
has attracted considerable attention. Classical questions include
quantum unique ergodicity, sup norm bounds, and the growth of
$L^4$-norms; see e.g. Holowinsky--Soundararajan \cite{HS}, Rudnick \cite{Rudnick}, and Blomer--Khan--Young \cite{BlomerKhanYoung2013distribution}.
Recently, Huang--Lester--Wigman--Yesha \cite{HLWY}
studied the sup norms of random holomorphic cusp forms.
Here we study extremal $L^p$-norms in the full space of cusp forms.
Peak forms attain the optimal orders on compact sets and on the whole
modular surface. Poincar\'e series provide a second explicit family
attaining the global order, and their norms admit precise asymptotics
when the index is at most a sufficiently small fixed multiple of the weight.

Let $\mathbb H=\{x+iy:y>0\}$, let $\Gamma=\SL_2(\mathbb Z)$, and put
\[
  \mathcal D=\{z=x+iy\in\mathbb H:|x|\leq 1/2,\ |z|\geq1\},
  \qquad \dd\mu(z)=\frac{\dd x\,\dd y}{y^2}.
\]
For an even integer $k\geq16$, write $\mathcal S_k$ for the space of
holomorphic cusp forms of weight $k$ for $\Gamma$. Thus
$f(\gamma z)=j(\gamma,z)^k f(z)$, where
$j(\gamma,z)=cz+d$ for $\gamma=\left(\begin{smallmatrix}a&b\\c&d\end{smallmatrix}\right)$.
We use the Petersson inner product
\[
  \langle f,g\rangle=\int_{\mathcal D}f(z)\overline{g(z)}y^k\,\dd\mu(z).
\]
For a measurable set $\mathcal K\subseteq\mathcal D$, define
\[
  \|f\|_{L^p(\mathcal K)}
  =\left(\int_{\mathcal K}|y^{k/2}f(z)|^p\,\dd\mu(z)\right)^{1/p}
  \quad(2\leq p<\infty),
  \qquad
  \|f\|_{L^\infty(\mathcal K)}
  =\sup_{z\in\mathcal K}|y^{k/2}f(z)|.
\]
We abbreviate $\|f\|_{L^2(\mathcal D)}$ to $\|f\|_2$.
All asymptotic statements refer to $k\to\infty$ through even integers;
$p$ is fixed, and $1/p=0$ when $p=\infty$.

\begin{theorem}\label{thm:main}
Let $\mathcal K\subset\mathcal D$ be a fixed compact set of positive
hyperbolic area. For each fixed $2<p\leq\infty$,
\begin{align}
  \sup_{\substack{f\in\mathcal S_k\\\|f\|_2=1}}
  \|f\|_{L^p(\mathcal K)}
  &\asymp_{\mathcal K,p} k^{1/2-1/p},\label{eq:local-extremal}\\
  \sup_{\substack{f\in\mathcal S_k\\\|f\|_2=1}}
  \|f\|_{L^p(\mathcal D)}
  &\asymp_p k^{3/4-3/(2p)}.\label{eq:global-extremal}
\end{align}
\end{theorem}

The upper bounds follow from
\[
  \|f\|_{L^\infty(\mathcal K)}\ll_{\mathcal K} k^{1/2}\|f\|_2,
  \qquad
  \|f\|_{L^\infty(\mathcal D)}\ll k^{3/4}\|f\|_2,
\]
due to Rudnick \cite{Rudnick} and Steiner \cite{Steiner}, respectively,
and the interpolation inequality
\[
  \|f\|_{L^p(\mathcal K)}
  \leq \|f\|_{L^\infty(\mathcal K)}^{1-2/p}\|f\|_2^{2/p}.
\]
For $p=\infty$, the sharpness of these bounds was
already explained in \cite[\S 2.1]{HLWY}. For $2<p<\infty$, we prove
the matching lower bounds by constructing the following forms.

\subsection{Peak forms}
Let $N_k=\dim\mathcal S_k$ and let $\{f_j\}_{j=1}^{N_k}$ be an
orthonormal basis. Set $F_j(z)=y^{k/2}f_j(z)$. The weighted Bergman kernel
\[
  \Pi_k(z,w)=\sum_{j=1}^{N_k}F_j(z)\overline{F_j(w)}
\]
is independent of the basis. Whenever $\Pi_k(w,w)>0$, define
\begin{equation}\label{eq:peak-definition}
  \Phi_{k,w}(z)=\frac{\Pi_k(z,w)}{\sqrt{\Pi_k(w,w)}},
  \qquad \phi_{k,w}(z)=y^{-k/2}\Phi_{k,w}(z).
\end{equation}
Indeed, by the definitions of $F_j$ and $\Pi_k$,
\[
  \phi_{k,w}(z)
  =\frac{1}{\sqrt{\Pi_k(w,w)}}
  \sum_{j=1}^{N_k}\overline{F_j(w)}f_j(z).
\]
Hence $\phi_{k,w}\in\mathcal S_k$, and orthonormality gives
\[
  \|\phi_{k,w}\|_2^2
  =\frac{\sum_{j=1}^{N_k}|F_j(w)|^2}{\Pi_k(w,w)}=1.
\]
We call $\phi_{k,w}$ a peak form with pole at $w$.

\begin{theorem}\label{thm:peakforms}
For every compact set $\mathcal K\subset\mathcal D$ of positive area,
there is a point $w_0\in\mathcal K$, independent of $k$ and $p$, such that
for every fixed $2<p\leq\infty$,
\[
  \|\phi_{k,w_0}\|_{L^p(\mathcal K)}
  \asymp_{\mathcal K,p} k^{1/2-1/p}.
\]
For $w_1=i(k-1)/(4\pi)$, we have
\[
  \|\phi_{k,w_1}\|_{L^p(\mathcal D)}
  \asymp_p k^{3/4-3/(2p)}.
\]
Both peak forms are well defined for all sufficiently large even $k$.
\end{theorem}

This proves Theorem~\ref{thm:main}. The local argument
uses the Bergman kernel estimates of
Huang--Lester--Wigman--Yesha \cite{HLWY}.
For the global assertion we give a short off-diagonal
calculation near the cusp. Related constructions of
bounded forms appear in \cite{HHH}.

\subsection{Poincar\'e series}
Put
\[
  \Gamma_\infty=
  \left\{\pm\begin{pmatrix}1&n\\0&1\end{pmatrix}:n\in\mathbb Z\right\},
  \qquad e(z)=e^{2\pi iz}.
\]
For $m\geq1$, the Poincar\'e series
\[
  P_{k,m}(z)=\sum_{\gamma\in\Gamma_\infty\backslash\Gamma}
  j(\gamma,z)^{-k}e(m\gamma z)
\]
belongs to $\mathcal S_k$; see \cite{iwaniec1997topics}.
When it is nonzero, write $\widetilde P_{k,m}=P_{k,m}/\|P_{k,m}\|_2$.
We first record bounds that hold in a wide range of indices.

\begin{theorem}\label{thm:wide}
Let $0<\varepsilon<1$. Uniformly for $1\leq m\leq k^{2-\varepsilon}$,
the form $P_{k,m}$ is nonzero for sufficiently large $k$, and
\begin{equation}\label{eq:wide-upper}
  \|\widetilde P_{k,m}\|_{L^\infty(\mathcal D)}
  \ll_\varepsilon
  \frac{k^{3/4}}{m^{1/2}}
  +k^\varepsilon\bigl(m^{1/4}+k^{5/12}\bigr).
\end{equation}
In the same range,
\begin{equation}\label{eq:wide-lower}
  \|\widetilde P_{k,m}\|_{L^\infty(\mathcal D)}
  \geq \bigl(2^{-5/4}\pi^{-3/4}+o_\varepsilon(1)\bigr)
  \frac{k^{3/4}}{m^{1/2}}.
\end{equation}
\end{theorem}

The nonvanishing of $P_{k,m}$ in this range
was already established by Rankin
\cite[Theorem~1]{RankinVanishing};
see also \cite[Proposition~2.2]{Ki2}.
The lower bound comes from the $m$-th Fourier coefficient, as in
Xia's method \cite{Xia}. It exceeds the universal bound
$\|\widetilde P_{k,m}\|_\infty\geq\vol(\mathcal D)^{-1/2}=\sqrt{3/\pi}$
when
\[
  m\leq\bigl(2^{-5/2}3^{-1}\pi^{-1/2}-\delta\bigr)k^{3/2}
\]
for any fixed sufficiently small $\delta>0$ and large $k$.

Our more precise result concerns the range in which the maximum of
$y^{k/2}e^{-2\pi my}$ lies a fixed distance above $y=1$.

\begin{theorem}\label{thm:asymptotic}
Fix $0<c_0<(4\pi)^{-1}$ and $2\leq p\leq\infty$.
Uniformly for $1\leq m\leq c_0k$, as $k\to\infty$ through even integers,
\begin{equation}\label{eq:exact-Lp}
  \|\widetilde P_{k,m}\|_{L^p(\mathcal D)}
  =C_p \frac{k^{3/4-3/(2p)}}{m^{1/2-1/p}}
  \bigl(1+O_{p,c_0}(k^{-1})\bigr),
\end{equation}
where
\[
  C_p=2^{3/p-5/4}\pi^{3/(2p)-3/4}p^{-1/(2p)}
  \quad(2\leq p<\infty),
  \qquad C_\infty=2^{-5/4}\pi^{-3/4}.
\]
In particular, $C_2=1$.
\end{theorem}

For any fixed $m$, this gives an explicit family attaining the global
order in Theorem~\ref{thm:main}. The proof compares the defining
Poincar\'e sum with its identity term in $L^p$, with an exponentially
small relative error. The case $p=2$ supplies the normalization, and
Stirling's formula gives \eqref{eq:exact-Lp}. The sup norm is treated
separately.

Kimmel studied the zeros of Poincar\'e series with $m\sim\alpha k$
in \cite{KimmelZeros}, and proved mass equidistribution in the range
$m/k\to\infty$, $m\ll k^{3/2-\varepsilon}$ in \cite{Kimmel}.
The latter paper also supplies the Fourier coefficient estimates used
for Theorem~\ref{thm:wide}. Theorem~\ref{thm:asymptotic} instead
describes the global norms in the range $m\leq c_0k$.


\medskip

For comparison, Sogge \cite{Sogge} established classical $L^p$-norm bounds for Laplace eigenfunctions on compact Riemannian manifolds. In the arithmetic setting, Iwaniec--Sarnak \cite{IwaniecSarnak} proved a power-saving sup norm bound for Hecke--Maass cusp forms. More recently, Humphries--Khan \cite{HumphriesKhan} obtained improved $L^p$-norm bounds, while Ki \cite{KiMaass} established the essentially optimal $L^4$-norm bound.

\section{Peak forms}\label{sec:peakforms}

For $z=x+iy$ and $w=u+iv$, put
\[
  \ell_\gamma(z,w)=
  \frac{2i\sqrt{yv}}{j(\gamma,z)(\gamma z-\overline w)},
  \qquad R_k(z,w)=\sum_{\gamma\in\Gamma}\ell_\gamma(z,w)^k.
\]
The standard Bergman kernel identity is
\begin{equation}\label{eq:Bergman}
  \Pi_k(z,w)=\frac{k-1}{8\pi}R_k(z,w);
\end{equation}
see \cite[(3.7)]{HLWY}. We will also use the elementary estimate
\begin{equation}\label{eq:lattice}
  \sum_{n\in\mathbb Z}\frac{1}{((n+a)^2+B^2)^2}
  \ll B^{-4}+B^{-3}
  \qquad(a\in\mathbb R,\ B>0),
\end{equation}
with an absolute constant. Indeed, a nonnegative unimodal function
has its sum over a translate of $\mathbb Z$ bounded by twice its
maximum plus its integral. Here the integral is $\pi/(2B^3)$.

\subsection{A compact set}
Choose a Lebesgue density point $w_0=u_0+iv_0$ of $\mathcal K$ in the
interior of $\mathcal D$, away from the elliptic points
$e_1=e^{\pi i/3}$, $e_2=i$, and $e_3=e^{2\pi i/3}$.
Such a point exists because $\mathcal K$ has positive area, by the
Lebesgue density theorem (see
\cite[Chapter~3, Corollary~1.5]{SteinShakarchi}).
Let $c_*>0$ be the absolute constant in \cite[Theorem 3.3]{HLWY},
which we may take at most $1$, and set
\[
  \delta=\min\left\{\frac12,\frac14\min_{1\leq j\leq3}|w_0-e_j|\right\},
  \qquad
  \rho=\min\left\{\frac12\operatorname{dist}(w_0,\partial\mathcal D),
    \frac{c_*\delta}{2}\right\}.
\]
For $|z-w_0|\leq\rho$, we have $z\in\mathcal D$,
$|z-e_j|>\delta$ for $1\leq j\leq3$, and
$|z-w_0|\leq c_*\delta$. Thus \cite[Theorem 3.3]{HLWY} gives,
uniformly in this explicitly specified neighborhood,
\begin{equation}\label{eq:local-kernel}
  R_k(z,w_0)=2\left(\frac{2i\sqrt{yv_0}}{z-\overline{w_0}}\right)^k
  +O_{\mathcal K}(e^{-c_{\mathcal K}k}),
\end{equation}
for some $c_{\mathcal K}>0$, since $y\asymp_{\mathcal K}1$ there.
Moreover,
\[
  |z-\overline{w_0}|^2=(x-u_0)^2+(y+v_0)^2
  =|z-w_0|^2+4yv_0,
\]
so
\[
  \left|\frac{2i\sqrt{yv_0}}{z-\overline{w_0}}\right|^k
  =\left(\frac{4yv_0}{|z-w_0|^2+4yv_0}\right)^{k/2}
  =\left(1+\frac{|z-w_0|^2}{4yv_0}\right)^{-k/2}.
\]
For all sufficiently large $k$, the disk
$|z-w_0|\leq v_0/\sqrt{k}$ lies in $|z-w_0|\leq\rho$,
and the last expression is bounded above and below by positive constants.
Using \eqref{eq:Bergman} and \eqref{eq:local-kernel}, we therefore have
$|\Pi_k(z,w_0)|\asymp_{\mathcal K}k$ throughout this disk, while
$\Pi_k(w_0,w_0)\sim(k-1)/(4\pi)$. By the definition
\eqref{eq:peak-definition},
\[
  |\Phi_{k,w_0}(z)|\asymp_{\mathcal K}k^{1/2}
  \qquad (|z-w_0|\leq v_0/\sqrt{k}).
\]
Since $w_0$ is a Lebesgue density point,
\[
  \operatorname{area}\bigl(\mathcal K\cap\{|z-w_0|\leq v_0/\sqrt{k}\}\bigr)
  =\frac{\pi v_0^2}{k}(1+o(1)).
\]
On this disk, $y^{-2}=v_0^{-2}(1+O(k^{-1/2}))$; hence
\[
  \mu\bigl(\mathcal K\cap\{|z-w_0|\leq v_0/\sqrt{k}\}\bigr)
  \sim\frac{\pi}{k}.
\]
Integration therefore proves
$\|\phi_{k,w_0}\|_{L^p(\mathcal K)}\gg_{\mathcal K,p}k^{1/2-1/p}$
for finite $p$. For $p=\infty$, evaluate at $w_0$.
The upper bounds follow from those stated in the introduction.

\subsection{The cusp}
We need the following estimate only on a strip of height comparable to $k$.

\begin{lemma}\label{lem:cusp-kernel}
Let $v=(k-1)/(4\pi)$, $w=u+iv$, and $z=x+iy$, where
$|x|,|u|\leq1/2$ and $|y-v|\leq\sqrt{k}$. There is an absolute $c>0$
such that, uniformly in this range,
\begin{equation}\label{eq:cusp-kernel}
  R_k(z,w)=
  2\frac{(4\pi\sqrt{yv})^k}{\Gamma(k)}
  e^{-2\pi(y+v)}e(x-u)\bigl(1+O(e^{-ck})\bigr).
\end{equation}
In particular, $|R_k(z,w)|\asymp\sqrt{k}$ and
$R_k(w,w)\sim\sqrt{2k/\pi}$.
\end{lemma}

\begin{proof}
First bound the contribution of $\gamma\in\Gamma\setminus\Gamma_\infty$.
For each nonidentity left coset choose a representative with primitive
bottom row $(c,d)$, $c\geq1$, and write $\gamma z=x_\gamma+iy_\gamma$.
Then $y_\gamma=y/|cz+d|^2\leq1/y$. Its representatives in $\Gamma$
are $\pm T^n\gamma$, where $Tz=z+1$. Since
\[
  |\ell_{T^n\gamma}(z,w)|^2
  =\frac{4vy_\gamma}{(n+x_\gamma-u)^2+(v+y_\gamma)^2}
  \leq\frac{4}{vy},
\]
equation \eqref{eq:lattice}, with $B=v+y_\gamma\geq v$, gives
\[
  \sum_{n\in\mathbb Z}|\ell_{T^n\gamma}(z,w)|^4
  \ll v^2y_\gamma^2\bigl((v+y_\gamma)^{-4}
  +(v+y_\gamma)^{-3}\bigr)
  \ll\frac{y_\gamma^2}{v},
\]
since $v\asymp k$. Next, using
$y_\gamma=y/((cx+d)^2+c^2y^2)$, equation \eqref{eq:lattice}
with $B=cy$ gives
\[
  \sum_{d\in\mathbb Z}y_\gamma^2
  =y^2\sum_{d\in\mathbb Z}
  \frac1{((cx+d)^2+c^2y^2)^2}
  \ll\frac1{c^4y^2}+\frac1{c^3y}.
\]
Summing over $c\geq1$ (and dropping the primitivity condition
on $(c,d)$), we obtain
\begin{align*}
  \sum_{\gamma\in\Gamma\setminus\Gamma_\infty}|\ell_\gamma(z,w)|^4
  &\ll\frac1v\sum_{c\geq1}
  \left(\frac1{c^4y^2}+\frac1{c^3y}\right)\\
  &\ll\frac1v(y^{-2}+y^{-1})\ll\frac1{vy}.
\end{align*}
Thus, by the bound $|\ell_\gamma(z,w)|^2\leq4/(vy)$ above,
\begin{align*}
  \sum_{\gamma\in\Gamma\setminus\Gamma_\infty}|\ell_\gamma(z,w)|^k
  &=\sum_{\gamma\in\Gamma\setminus\Gamma_\infty}
  |\ell_\gamma(z,w)|^4
  \bigl(|\ell_\gamma(z,w)|^2\bigr)^{k/2-2}\\
  &\leq\left(\frac4{vy}\right)^{k/2-2}
  \sum_{\gamma\in\Gamma\setminus\Gamma_\infty}|\ell_\gamma(z,w)|^4\\
  &\ll\frac1{vy}\left(\frac4{vy}\right)^{k/2-2}
  \ll (C/k)^{k-2},
\end{align*}
for an absolute $C>0$, since $y,v\asymp k$.

The parabolic terms give
\[
  R_k(z,w)=2(2i\sqrt{yv})^k
  \sum_{n\in\mathbb Z}(z-\overline w+n)^{-k}
  +O((C/k)^{k-2}).
\]
Using the classical identity
\[
  \sum_{n\in\mathbb Z}(z+n)^{-k}
  =\frac{(-2\pi i)^k}{\Gamma(k)}\sum_{r\geq1}r^{k-1}e(rz)
  \qquad(z\in\mathbb H),
\]
see \cite[p.~92]{serre1973course}, we obtain
\[
  R_k(z,w)=2\frac{(4\pi\sqrt{yv})^k}{\Gamma(k)}
  \sum_{r\geq1}r^{k-1}e^{-2\pi(y+v)r}e(r(x-u))
  +O((C/k)^{k-2}).
\]
Here $2\pi(y+v)/(k-1)=1+O(k^{-1/2})$. For large $k$ this ratio
is at least $9/10$, whereas $\log r\leq(r-1)\log2$ for $r\geq2$.
Consequently, for $r\geq2$,
\[
  r^{k-1}e^{-2\pi(y+v)(r-1)}
  \leq e^{-(9/10-\log2)(k-1)(r-1)}.
\]
Since $9/10-\log2>0$, summing the geometric series gives
\[
  \sum_{r\geq2}r^{k-1}e^{-2\pi(y+v)(r-1)}\ll e^{-c k}.
\]
The $r=1$ term therefore dominates.

To bound its size, put
$A(y)=(4\pi\sqrt{yv})^ke^{-2\pi(y+v)}/\Gamma(k)$.
Stirling's formula gives $A(v)\sim\sqrt{k/(2\pi)}$. If
$h=y-v$ and $|h|\leq\sqrt{k}$, then
\[
  \log\frac{A(v+h)}{A(v)}
  =\frac{k}{2}\log(1+h/v)-2\pi h
  =\frac{h}{2v}+O\left(\frac{kh^2}{v^2}\right)=O(1).
\]
Thus $A(y)\asymp\sqrt{k}$; the nonparabolic error can be absorbed
into the relative exponential error in \eqref{eq:cusp-kernel}.
\end{proof}

Take $w=w_1=iv$. Equations \eqref{eq:Bergman} and
\eqref{eq:cusp-kernel} show that
$\Pi_k(w,w)\asymp k^{3/2}$ and
$|\Phi_{k,w}(z)|\asymp k^{3/4}$ on
$|x|\leq1/2$, $|y-v|\leq\sqrt{k}$. This strip lies in $\mathcal D$
for large $k$ and has hyperbolic area $\asymp k^{-3/2}$. Hence
\[
  \|\phi_{k,w}\|_{L^p(\mathcal D)}^p
  \gg_p k^{3p/4-3/2}\qquad(2<p<\infty).
\]
For $p=\infty$, evaluate at $w$.
The general upper bounds complete the proof of Theorem~\ref{thm:peakforms}.

\section{Exact norms of Poincar\'e series}\label{sec:exact}

We prove Theorem~\ref{thm:asymptotic} directly from the defining sum.
Throughout this section, $0<c_0<(4\pi)^{-1}$ is fixed and
$1\leq m\leq c_0k$. For $z=x+iy$ and $\gamma\in\Gamma$, write
$y_\gamma=\Im(\gamma z)$. Put
\[
  H(y)=y^{k/2}e^{-2\pi my},\qquad
  y_0=\frac{k}{4\pi m},\qquad
  s=\frac{1-4\pi c_0}{2}\,k.
\]
We take $k$ sufficiently large that $s\geq2$.
The function $H$ is maximal at $y_0>1$.
For $0<y\leq Y\leq1$, we have
\[
  \log\frac{H(Y)}{H(y)}
  =\int_y^Y\left(\frac{k}{2}-2\pi m u\right)
    \frac{\dd u}{u}
  \geq s\log\frac{Y}{y},
\]
since $u(\log H)'(u)=k/2-2\pi mu\geq s$ on $(0,1]$. Hence
\begin{equation}\label{eq:H-comparison}
  H(y)\leq H(Y)(y/Y)^s.
\end{equation}

Also, $y_0\geq(4\pi c_0)^{-1}>1$, and the function
$\log u+u^{-1}-1$ is increasing for $u>1$. Thus, with
$\kappa=\tfrac12(\log(1/(4\pi c_0))+4\pi c_0-1)>0$,
\begin{equation}\label{eq:exponential-gap}
  \frac{H(1)}{H(y_0)}
  =\exp\left\{\frac{k}{2}(1-\log y_0-y_0^{-1})\right\}
  \leq e^{-\kappa k}.
\end{equation}

\subsection{The nonidentity terms}
Write $\mathscr C=\Gamma_\infty\backslash\Gamma$,
$\mathbf1=\Gamma_\infty I$, and
$\mathscr C^*=\mathscr C\setminus\{\mathbf1\}$.
The value $y_\gamma$ depends only on the coset. On $\mathcal D$
we have $y_\gamma\leq y$,
since $|cz+d|\geq1$ for every primitive bottom row $(c,d)$.
For $\gamma\in\mathscr C^*$, choose $c\geq1$; then also
$y_\gamma\leq1/y$.

We require only fixed fourth-power lattice sums. From
\eqref{eq:lattice}, enlarging the sum over primitive rows to all rows,
\begin{equation}\label{eq:height-square}
  \sum_{\gamma\in\mathscr C^*}y_\gamma^2
  \leq\sum_{c\geq1}\sum_{d\in\mathbb Z}
  \frac{y^2}{((cx+d)^2+c^2y^2)^2}
  \ll y^{-2}+y^{-1}.
\end{equation}
Consequently, for every real $s\geq2$ (in particular for the value
fixed above), we have, with absolute constants,
\begin{align}
  \sum_{\gamma\in\mathscr C^*}y_\gamma^s&\ll1
  &&(z\in\mathcal D,\ y\leq1),\label{eq:compact-sum}\\
  \sum_{\gamma\in\mathscr C^*}(yy_\gamma)^s&\ll y
  &&(z\in\mathcal D,\ y\geq1).\label{eq:cusp-sum}
\end{align}
Indeed, in the first range $0<y_\gamma\leq1$, and in the second
$0<yy_\gamma\leq1$, so their $s$-th powers are bounded by their
squares. In the first range $y\geq\sqrt3/2$; in the second, multiply
\eqref{eq:height-square} by $y^2$.

Decompose the weighted sum as
\[
  y^{k/2}P_{k,m}(z)=E(z)+R(z),\qquad E(z)=y^{k/2}e(mz).
\]
The weighted $\gamma$-summand has absolute value
\[
  y^{k/2}|j(\gamma,z)|^{-k}e^{-2\pi my_\gamma}=H(y_\gamma).
\]
Applying \eqref{eq:H-comparison} with $Y=1$ or $Y=1/y$, and then
\eqref{eq:compact-sum}--\eqref{eq:cusp-sum}, gives
\begin{equation}\label{eq:remainder-pointwise}
  |R(z)|\ll
  \begin{cases}
  H(1),& y\leq1,\\
  yH(1/y)\leq H(1)y^{1-s},& y\geq1.
  \end{cases}
\end{equation}
Since $\mu(\mathcal D)<\infty$ and $s\geq2$, we obtain
\begin{equation}\label{eq:remainder-norms}
  \left(\int_{\mathcal D}|R(z)|^p\,\dd\mu(z)\right)^{1/p}
  \ll_p H(1)\quad(2\leq p<\infty),
  \qquad \sup_{\mathcal D}|R(z)|\ll H(1).
\end{equation}

\subsection{Finite moments and normalization}
For fixed $2\leq p<\infty$, define the full gamma integral
\begin{equation}\label{eq:gamma-integral}
  A_p=\int_0^\infty H(y)^p\frac{\dd y}{y^2}
  =\frac{\Gamma(kp/2-1)}{(2\pi mp)^{kp/2-1}}.
\end{equation}
The fundamental domain is contained in the vertical strip
$|x|\leq1/2$ and contains its part with $y\geq1$. Therefore
\begin{equation}\label{eq:seed-integral}
  0\leq A_p-\int_{\mathcal D}|E(z)|^p\,\dd\mu(z)
  \leq\int_0^1H(y)^p\frac{\dd y}{y^2}
  \leq\frac{H(1)^p}{ps-1}.
\end{equation}
The last inequality uses $H(y)\leq H(1)y^s$ and explicitly accounts
for the factor $y^{-2}$ near zero.

We next compare $H(1)$ with $A_p^{1/p}$.
For $0\leq u\leq k^{-1/2}$,
\[
  \log\frac{H(y_0(1+u))}{H(y_0)}
  =\frac{k}{2}\bigl(\log(1+u)-u\bigr)\geq-\frac14.
\]
Integrating over this interval and using $y_0\leq k/(4\pi)$ gives
\begin{equation}\label{eq:saddle-lower}
\begin{aligned}
  A_p&\geq\int_{y_0}^{y_0(1+k^{-1/2})}
  H(y)^p\frac{\dd y}{y^2}
  \gg_p\frac{H(y_0)^p}{y_0\sqrt{k}}
  \gg_p k^{-3/2}H(y_0)^p.
\end{aligned}
\end{equation}
By \eqref{eq:exponential-gap} and \eqref{eq:saddle-lower},
\begin{equation}\label{eq:relative-gap}
  \frac{H(1)}{A_p^{1/p}}\ll_p k^{3/(2p)}e^{-\kappa k}.
\end{equation}
Combining \eqref{eq:seed-integral} with \eqref{eq:relative-gap}, we obtain
\[
  \int_{\mathcal D}|E(z)|^p\,\dd\mu(z)
  =A_p+O\left(\frac{H(1)^p}{ps-1}\right)
  =A_p\bigl(1+O_{p,c_0}(k^{1/2}e^{-p\kappa k})\bigr).
\]
By Minkowski's inequality, \eqref{eq:remainder-norms}, and
\eqref{eq:relative-gap},
\[
  \left|\|P_{k,m}\|_{L^p(\mathcal D)}
  -\left(\int_{\mathcal D}|E(z)|^p\,\dd\mu(z)\right)^{1/p}\right|
  \ll_p H(1)\ll_p k^{3/(2p)}e^{-\kappa k}A_p^{1/p}.
\]
It follows that, for some $c=c(c_0)>0$,
\begin{equation}\label{eq:unnormalized-norm}
  \|P_{k,m}\|_{L^p(\mathcal D)}
  =A_p^{1/p}\bigl(1+O_{p,c_0}(e^{-ck})\bigr).
\end{equation}
In particular, taking $p=2$ proves both nonvanishing and
\begin{equation}\label{eq:linear-normalization}
  \|P_{k,m}\|_2^2
  =\frac{\Gamma(k-1)}{(4\pi m)^{k-1}}
  \bigl(1+O_{c_0}(e^{-ck})\bigr).
\end{equation}
No estimate of Fourier coefficients is needed for this normalization.

Dividing \eqref{eq:unnormalized-norm} by its $p=2$ case, and using
\eqref{eq:gamma-integral}, yields
\begin{equation}\label{eq:normalized-gamma}
  \|\widetilde P_{k,m}\|_{L^p(\mathcal D)}^p
  =\frac{(4\pi m)^{p(k-1)/2}}{\Gamma(k-1)^{p/2}}
  \frac{\Gamma(kp/2-1)}{(2\pi mp)^{kp/2-1}}
  \bigl(1+O_{p,c_0}(e^{-ck})\bigr).
\end{equation}
For fixed $a>0$, Stirling's formula in the form
\[
  \Gamma(ak-1)=\sqrt{2\pi}(ak)^{ak-3/2}e^{-ak}
  \bigl(1+O_a(k^{-1})\bigr)
\]
reduces the right side of \eqref{eq:normalized-gamma} to
\[
  2^{3-5p/4}\pi^{3/2-3p/4}p^{-1/2}
  k^{3p/4-3/2}m^{1-p/2}
  \bigl(1+O_{p,c_0}(k^{-1})\bigr).
\]
Taking the $p$-th root proves \eqref{eq:exact-Lp} for finite $p$.

\subsection{The sup norm}
Since $y_0>1$, the identity term satisfies
$\sup_{z\in\mathcal D}|E(z)|=H(y_0)$, attained at $z=iy_0$.
Here it is the ratio $H(1)/H(y_0)$, not $H(1)$ itself, that
is exponentially small. Indeed, \eqref{eq:remainder-norms} and
\eqref{eq:exponential-gap} give, independently of the finite-$p$ argument,
\[
  \left|\|P_{k,m}\|_{L^\infty(\mathcal D)}-H(y_0)\right|
  \leq\sup_{\mathcal D}|R|\ll H(1)
  \ll H(y_0)e^{-\kappa k}.
\]
Using \eqref{eq:linear-normalization}, we obtain
\begin{align*}
  \|\widetilde P_{k,m}\|_{L^\infty(\mathcal D)}
  &=\frac{(4\pi m)^{(k-1)/2}}{\Gamma(k-1)^{1/2}}
    \left(\frac{k}{4\pi m}\right)^{k/2}e^{-k/2}
    \bigl(1+O_{c_0}(e^{-ck})\bigr)\\
  &=2^{-5/4}\pi^{-3/4}\frac{k^{3/4}}{m^{1/2}}
    \bigl(1+O_{c_0}(k^{-1})\bigr),
\end{align*}
by Stirling's formula. This completes the proof of
Theorem~\ref{thm:asymptotic}.

\section{Bounds for larger indices}\label{sec:wide}

\subsection{Fourier coefficients}
The Fourier expansion is
\begin{equation}\label{eq:Fourier}
  P_{k,m}(z)=\sum_{n\geq1}p_{k,m}(n)e(nz),
\end{equation}
where
\[
  p_{k,m}(n)=\delta_{m,n}+2\pi i^k
  \left(\frac nm\right)^{(k-1)/2}
  \sum_{c\geq1}\frac{S(m,n;c)}c
  J_{k-1}\left(\frac{4\pi\sqrt{mn}}c\right).
\]
Here $S(m,n;c)=\sideset{}{^{*}}\sum_{d\bmod c}e((md+n\overline d)/c)$
is the Kloosterman sum, and $J_{k-1}$ is the Bessel function.
We use the following consequence of Kimmel's estimates.

\begin{lemma}\label{lem:coefficients}
For even $k\geq16$, $m,n\geq1$, and every $\eta>0$,
\begin{equation}\label{eq:coefficients}
  |p_{k,m}(n)-\delta_{m,n}|
  \ll_\eta\left(\frac nm\right)^{(k-1)/2}
  \left\{\frac{(mn)^{1/4+\eta}}k
  +\frac{(mn)^\eta}{k^{1/3}}+\frac{\sqrt{mn}}{k^2}\right\}.
\end{equation}
In particular, for every fixed $0<\varepsilon<1$,
$p_{k,m}(m)=1+o_\varepsilon(1)$ uniformly for $m\leq k^{2-\varepsilon}$.
\end{lemma}

\begin{proof}
Put $Q=4\pi\sqrt{mn}/(k-1)$ and $g=(m,n)$.
For $c\leq Q$, \cite[Lemma 2.5]{Kimmel} bounds the Bessel--Kloosterman
sum by
\[
  \ll_\eta (mn)^{\eta/2}
  \left\{\frac{(mn)^{1/4}}k\sum_{d\mid g}d^{-1/2}
  +\frac{k^{1/6}}{(mn)^{1/4}}
  \sum_{\substack{d\mid g\\d\leq Q}}d^{1/2}\right\}.
\]
Use $\sum_{d\mid g}d^{-1/2}\leq\tau(g)\ll_\eta(mn)^{\eta/2}$ and
$\sum_{d\mid g,\ d\leq Q}d^{1/2}\leq Q^{1/2}\tau(g)$.
For $c>Q$, \cite[Lemma 2.6, (2.8)]{Kimmel} gives
$\ll Q/k+k^{-1/3}$. These estimates prove \eqref{eq:coefficients}.
Putting $n=m$ and choosing
$0<\eta<\min(\varepsilon/16,1/48)$ proves the last assertion.
\end{proof}

Unfolding the Petersson inner product gives the exact identity
\begin{equation}\label{eq:norm-exact}
  \|P_{k,m}\|_2^2
  =\frac{\Gamma(k-1)}{(4\pi m)^{k-1}}p_{k,m}(m).
\end{equation}
Thus the last part of Lemma~\ref{lem:coefficients} gives nonvanishing
and the normalization throughout the range of Theorem~\ref{thm:wide}.

\subsection{The upper bound}
We use a single summation estimate to handle the Fourier series.
For fixed $0\leq a\leq1/2$, put
\[
  T_a(y)=\frac{y^{k/2}}{(k/e)^{k/2}}
  \sum_{n\geq1}(4\pi n)^{k/2-a}e^{-2\pi ny}.
\]
Comparison of a unimodal sum with its integral and maximum, followed
by Stirling's formula, gives
\begin{equation}\label{eq:Ta}
  T_a(y)\ll
  k^{1/2-a}y^{a-1}+k^{-a}y^a.
\end{equation}
Indeed, before multiplying by $y^{k/2}/(k/e)^{k/2}$ the integral is
\[
  \frac{(4\pi)^{k/2-a}\Gamma(k/2-a+1)}{(2\pi y)^{k/2-a+1}},
\]
and the maximum is
$((k-2a)/y)^{k/2-a}e^{-k/2+a}$.
For $y\geq k$, there is also the stronger bound
\begin{equation}\label{eq:Ta-tail}
  T_a(y)\ll e^{-4k}.
\end{equation}
To see this, use $\log n\leq n-1$ to bound the sum by a constant
times $(4\pi)^{k/2-a}e^{-2\pi y}$. The logarithm of the remaining
factor is $\tfrac k2\log(4\pi e y/k)-2\pi y$, which decreases for
$y\geq k$ and is less than $-4k$ at $y=k$.

Assume $m\leq k^{2-\varepsilon}$, and fix
$0<\eta<\min(\varepsilon/16,1/48)$.
By Lemma~\ref{lem:coefficients}, \eqref{eq:norm-exact},
and Stirling's formula,
\[
  \|P_{k,m}\|_2^{-1}\asymp
  \frac{(4\pi m)^{(k-1)/2}}{\Gamma(k-1)^{1/2}},
  \qquad
  \frac{(k/e)^{k/2}}{\Gamma(k-1)^{1/2}}\asymp k^{3/4}.
\]
Consequently, for $0\leq a\leq1/2$,
\[
  \frac{y^{k/2}}{\|P_{k,m}\|_2}
  \sum_{n\geq1}\left(\frac nm\right)^{(k-1)/2}
  (mn)^{1/2-a}e^{-2\pi ny}
  \ll m^{1/2-a}k^{3/4}T_a(y).
\]
Taking $a=1/4-\eta,1/2-\eta,0$, with prefactors
$k^{-1},k^{-1/3},k^{-2}$, handles \eqref{eq:coefficients}.
The $\delta_{m,n}$ term is $\ll k^{3/4}m^{-1/2}$,
since $y^{k/2}e^{-2\pi my}$ is maximal at $y=k/(4\pi m)$.
Thus \eqref{eq:Fourier} gives
\begin{align}\label{eq:upper-components}
  |y^{k/2}\widetilde P_{k,m}(z)|
  \ll_\varepsilon{}&\frac{k^{3/4}}{m^{1/2}}
  +m^{1/4+\eta}k^{-1/4}T_{1/4-\eta}(y)\nonumber\\
  &+m^\eta k^{5/12}T_{1/2-\eta}(y)
  +m^{1/2}k^{-5/4}T_0(y).
\end{align}
For $\sqrt3/2\leq y\leq k$, \eqref{eq:Ta} gives respectively
\begin{align*}
  m^{1/4+\eta}k^{-1/4}T_{1/4-\eta}(y)
  &\ll m^{1/4}k^{3\eta},\\
  m^\eta k^{5/12}T_{1/2-\eta}(y)
  &\ll k^{5/12+3\eta},\\
  m^{1/2}k^{-5/4}T_0(y)
  &\ll k^{1/4}.
\end{align*}
For $y\geq k$, all three terms are exponentially small by
\eqref{eq:Ta-tail}. Since $3\eta<\varepsilon$, this proves
\eqref{eq:wide-upper}.

\subsection{The lower bound}
For any $y>0$, Fourier coefficient extraction gives
\[
  I(y):=\int_{-1/2}^{1/2}y^{k/2}\widetilde P_{k,m}(x+iy)e(-mx)\,\dd x
  =\frac{p_{k,m}(m)}{\|P_{k,m}\|_2}y^{k/2}e^{-2\pi my}.
\]
The weighted absolute value of a modular form is $\Gamma$-invariant.
Hence $|I(y)|\leq\|\widetilde P_{k,m}\|_{L^\infty(\mathcal D)}$
even when the horocycle is not contained in $\mathcal D$.
At $y_m=k/(4\pi m)$, \eqref{eq:norm-exact} gives
\begin{align*}
  I(y_m)
  &=p_{k,m}(m)^{1/2}
  \frac{(4\pi m)^{(k-1)/2}}{\Gamma(k-1)^{1/2}}
  \left(\frac{k}{4\pi m}\right)^{k/2}e^{-k/2}\\
  &=\bigl(2^{-5/4}\pi^{-3/4}+o_\varepsilon(1)\bigr)
  \frac{k^{3/4}}{m^{1/2}}.
\end{align*}
Here $p_{k,m}(m)>0$ by \eqref{eq:norm-exact}, and the last step uses
Lemma~\ref{lem:coefficients} and Stirling's formula.
This proves \eqref{eq:wide-lower} and Theorem~\ref{thm:wide}.

\section*{Acknowledgements}
The study of the sup norm of Poincar\'e series in this paper was suggested by Ze\'ev Rudnick. We thank him for this suggestion and his interest in this work. All mathematical arguments in this paper were completed when the first-named author (X. Han) visited the second-named author (B. Huang) at Shandong University in June 2026. X. Han would like to thank Shandong University for its warm hospitality.

\medskip
\noindent
\textbf{Use of AI.} ChatGPT was used to improve the language of the manuscript. The authors take full responsibility for the mathematical results and proofs.

\end{document}